\documentclass[preprint,12pt]{elsarticle}
\usepackage{color}

\usepackage{amsmath}
\usepackage{amsthm}
\usepackage{amsfonts}
\usepackage{amssymb}
\usepackage{graphicx}
\usepackage{eepic}
\usepackage{epsfig} 
\usepackage{amsbsy}
\usepackage{enumerate}
\usepackage{verbatim}
\usepackage{epstopdf}
\usepackage{subfigure}
\graphicspath{{../figure/}}
\newtheorem{remark}{Remark}[section]

\newcommand{\E}{{\mathbb E}}
\newcommand{\cO}{{\mathcal O}}

\newcommand{\HH}{{\mathcal H}}

\newcommand{\cN}{{\mathcal N}}

\newcommand{\cP}{{\mathcal P}}
\newcommand{\cH}{{\mathcal H}}

\newcommand{\R}{{\mathbb R  }}

\newcommand{\cA}{\mathcal A}
\newcommand{\cL}{\mathcal L}

\newcommand{\eps}{\epsilon}

\journal{Journal of Computational Physics}

\begin{document}

\begin{frontmatter}

\title{Numerical methods for stochastic partial differential equations with multiples scales}
%\author{A. Abdulle\footnote{Mathematics Section, Ecole Polytechnique F\'ed\'erale de Lausanne, CH-1015 Lausanne, Switzerland}~ and G.A. Pavliotis\footnote{Department of Mathematics, Imperial College London, London SW7 2AZ, UK}
%                    }

\author{A. Abdulle}
\address{Mathematics Section, Ecole Polytechnique F\'ed\'erale de Lausanne, CH-1015 Lausanne, Switzerland}

\author{G.A. Pavliotis}
\address{Department of Mathematics, Imperial College London, London SW7 2AZ, UK}

\begin{abstract}
A new method for  solving numerically stochastic partial differential equations (SPDEs) with multiple scales is presented.
The method combines a spectral method with the heterogeneous multiscale method (HMM) presented in [W.~E, D.~Liu, and E.~Vanden-Eijnden, Comm. Pure Appl. Math., 58(11):1544--1585, 2005]. The class of problems that we consider are  SPDEs with quadratic nonlinearities that were studied in [D.~Bl{\"o}mker, M.~Hairer, and G.~A. Pavliotis,  Nonlinearity, 20(7):1721--1744, 2007.]
For such SPDEs an amplitude equation which describes the effective dynamics at long time scales can be rigorously derived for both 
advective and diffusive time scales. Our method, based on micro and macro solvers, allows to capture numerically the amplitude equation accurately at a cost
independent of the small scales in the problem. Numerical experiments illustrate the behavior of the proposed method.
\end{abstract}  


\end{frontmatter}
%
%%%%%%%%%%%%%%%%%%%%%%%%%%%%%%%%%%%%%%%%%%%%%%%%%%%%%%%%%%%%%%%%%%%%%%%%%%%%%%%
%
\section{Introduction}~\label{sec:intro}
Many interesting phenomena in the physical sciences and in applications are characterized by their high dimensionality and the presence of many different spatial and temporal scales.  Standard examples include atmosphere and ocean sciences~\cite{MajdaFranzkeKhouider2008}, molecular dynamics~\cite{Griebel} and materials science~\cite{Fish}. 
The mathematical description of phenomena of this type quite often leads to infinite dimensional multiscale systems that are described by nonlinear evolution partial differential equations (PDEs) with multiple scales.

Often physical systems are also subject to noise. This noise might be either due to thermal fluctuations~\cite{Einst56}, noise in some control parameter~\cite{HorsLef84}, coarse-graining of a high-dimensional deterministic system with random initial conditions~\cite{Mazo02, Zwan01}, or the stochastic parameterization of small scales~\cite{ELV05}. High dimensional multiscale dynamical systems that are subject to noise can be modeled accurately using stochastic partial differential equations (SPDEs) with a multiscale structure. There are very few instances where 
SPDEs with multiple scales can be treated analytically. The goal of this paper is to develop numerical methods for solving accurately and efficiently multiscale SPDEs. Several numerical methods for SPDEs have been developed and analyzed in recent  years, 
e.g.~\cite{AlabertGyongy06, DavieGaines01, Printems01}, based on a finite difference scheme in both space and time.
% A finite difference scheme in space has been studied in \cite{AlabertGyongy06} for the one dimensional Burgers equation. A finite difference scheme in space and time 
%has been studied in \cite{DavieGaines01} for one dimensional parabolic problems 
%and a finite difference scheme in space with a semi-implicit method in time has
%been studied in \cite{Printems01}. 
It is well known that explicit time discretization via standard methods (e.g., as the Euler-Maruyama method) leads to a time-step restriction due to the stiffness originating from the discretisation of the diffusion operator (e.g. the CFL condition $\Delta t\leq C(\Delta x)^2$, where
$\Delta t$ and $\Delta x$ are the time and space discretization, respectively).
The situation is even worse for SPDEs with multiple scales (e.g. of the form~\eqref{e:spde_intro_adv} and~\eqref{e:spde_intro_diff} below) as in this case the Courant-Friedrichs-Lewy (CFL) condition becomes $\Delta t\leq C(\Delta x\cdot\eps)^2$, where $\eps \ll 1$ is the parameter measuring scale separation. Standard
explicit methods become useless for SPDEs with multiple scales. 

Such time-step restriction can in theory be removed by using implicit methods 
as was shown in \cite{Printems01}. However the implicitness of the numerical scheme forces one to solve  potentially large linear algebraic problems at each time step. Furthermore, it was shown in~\cite{Li08} that implicit methods are not suited for studying the long time dynamics of fast-slow stochastic systems as they do not capture the correct invariant measure of the system. Although this result has been obtained for finite dimensional stochastic systems, it is expected that it also applies to infinite dimensional fast-slow systems of stochastic differential equations (SDEs), rendering the applicability of implicit methods to SPDEs with multiple scales questionable.
We also note that a new class of explicit methods, the S-ROCK methods, with much better stability properties than the Euler-Maruyama method was recently introduced in~\cite{Abd07, Abd08, Abd08b}. Although these methods are much more efficient than traditional explicit methods, computing time issues will occur when trying to solve SPDEs with multiple scales as considered here, since the stiffness is extremely severe for small $\eps.$ Furthermore, capturing the correct invariant measure of the SPDE for $\Delta t>\eps$ is still an issue for such solvers.

In this paper we consider SPDEs of the form
\begin{equation}
\label{e:spde}
 \partial_t v = \cA v + F(v) + \eps \, Q\xi,
\end{equation}
posed in a bounded domain of $\R$ with appropriate boundary conditions. The differential operator $\cA$
 is assumed to be a non-positive self-adjoint operator in a Hilbert space $\cH$, $\xi$ denotes space-time Gaussian white noise, $Q$ is the covariance operator of the noise and we take $\eps \ll 1$. We assume that the operator $\cA$ has a finite dimensional kernel, $\cN$. This assumption leads to scale separation between the slow dynamics in $\cN$ and the fast dynamics in the orthogonal complement of the null space $\cN^{\bot},$ where $\cH = \cN \oplus \cN^{\bot}$. In this paper we will furthermore assume that noise acts directly only on the orthogonal complement $\cN^{\bot}$. When noise acts also on $\cN$, different distinguished limits than the ones considered in this paper should be considered.

In order to describe the longtime behavior of the SPDEs 
We perform an \emph{advective} rescaling set $v(t):=\eps u(\eps t)$. Using the scaling properties of white noise we obtain the following singularly perturbed SPDE
\begin{equation}\label{e:spde_intro_adv}
\partial_t u = \frac{1}{\eps} \cA u + F(u)+ \frac{1}{\sqrt{\eps}} Q\xi. 
\end{equation}
Another scaling is of interest, namely the \emph{diffusive} rescaling $v(t):=\eps u(\eps^2 t)$ 
%and $F(v)=f(v)+\eps^2 g(v)$,
which leads to the SPDE
\begin{equation}\label{e:spde_intro_diff}
\partial_t u = \frac{1}{\eps^2} \cA u + \frac{1}{\eps} F(u)+ \frac{1}{\eps} Q\xi. 
\end{equation}
For concreteness, we will focus on the class of SPDEs with quadratic nonlinearities  that was considered in~\cite{BHP06}, and assume that $F(u)=f(u)+\eps^{\alpha} g(u),~$
where $f$ is a quadratic function (e.g. $f(u)=B(u,u)$, a symmetric bilinear form), $g$ a linear function and the exponent $\alpha$ is either $1$ or $2$.\footnote{Usually the functions $f$ and $g$ involve derivatives of the function $u$. For example, for both the Burgers and the Kuramoto-Shivashinsky equation we have $f(u) = u \partial_x u$. The linear function $g(u)$ is included to induce a linear instability to the dynamics. In the case of the Burgers equation we will simply take $g(u) = u$ whereas in the case of the Kuramoto-Shivashinsky equation we can take $g(u) = \partial_x^4 u$. Further discussion can be found in Section~\ref{sec:numerics} and in~\cite{Pavliotis_al2010}. }  The choice of $\alpha$ will depend on the particular scaling. 
Singularly perturbed SPDEs with quadratic nonlinearities provide a natural testbed for testing the applicability of the heterogeneous multiscale method to infinite dimensional stochastic systems, since a rigorous homogenization theory exists for this class of SPDEs~\cite{BHP06}. Furthermore, SPDEs of this form arise naturally in stochastic models for climate~\cite{MajdaFranzkeKhouider2008} and in surface growth~\cite{KPZ1986, zaleski87}. Finally, it has already been shown through rigorous analysis and numerical experiments that these systems exhibit a very rich dynamical behavior, such as noise-induced transitions~\cite{WEal2010} and the possibility of stabilization of linearly unstable modes due to the interaction between the additive noise and the scale separation~\cite{DBMHGP08, Pavliotis_al2010}. We believe, however, that the methodology developed in this paper has a wider range of applicability and is not restricted to SPDEs with quadratic nonlinearities. Further comments about the class of SPDEs for which we believe that the proposed numerical method can be applied can be found in Section~\ref{sec:conclusions}.

Our numerical algorithm is based on a combination of a spectral method with micro-macro time integration schemes.
We denote by $x=P_c u$ the projection onto $\cN$  and by $y=P_s u$, $P_s=I-P_c$ the projection onto $\cN^{\bot}$. 
We then rewrite~\eqref{e:spde_intro_adv} and~\eqref{e:spde_intro_diff} as fast-slow system of SDEs
\begin{subequations}\label{e:fast_slow_intro_adv}
\begin{eqnarray}
\dot{x} & = &  a(x,y), \\
\dot{y}  & = & \frac{1}{\eps} \cA y + b(x,y) + \frac{1}{\sqrt{\eps}} Q\xi,
\end{eqnarray}
\end{subequations}
and
\begin{subequations}\label{e:fast_slow_intro_diff}
\begin{eqnarray}
\dot{x} & = & \frac{1}{\eps} a(x,y), \label{e:xx} \\
\dot{y}  & = & \frac{1}{\eps^2} \cA y + \frac{1}{\eps} b(x,y) + \frac{1}{\eps}
Q\xi,
\end{eqnarray}
\end{subequations}
where the functions $a(x,y)$ and  $b(x,y)$ are the projections of $F(u)$ onto $\cN$ and $\cN^{\bot}$. We remark that an ${\mathcal{O}(1)}$ nonlinear term can be added in~\eqref{e:xx}. The fast-slow systems~\eqref{e:fast_slow_intro_adv} and~\eqref{e:fast_slow_intro_diff} resemble fast slow systems for SDEs \cite[Ch. 10,11]{PavlSt08}.
%to which standard homogenization and averaging can be applied \cite[Ch. 10,11]{PavlSt08}.
However, 
%unlike equations~\eqref{e:averaging} and~\eqref{e:homogenization},
 the fast process $y$ is infinite dimensional and the well known averaging and homogenization theorems~\cite{lions, PapStrVar77}
do not apply.

Averaging and homogenization results for SPDEs have been obtained recently~\cite{CerraiFreidlin09, BHP06}. In particular, provided that the fast process $y$ in~\eqref{e:fast_slow_intro_adv} has suitable ergodic properties, then the slow variable $x$ converges, in the limit as $\eps$ tends to $0,$ to the solution of the averaged equation
\begin{equation}\label{e:aver}
\dot{x} = \bar{a}(x),
\end{equation}
where the averaged coefficient is given by the average of $a(x, y)$ with respect to the invariant measure of the (infinite dimensional) fast process $y$. When this average vanishes (i.e. the centering condition from homogenization theory is satisfied) then the dynamics at the advective time scale becomes trivial and it is necessary to look at the dynamics at the diffusive time scale, equations~\eqref{e:fast_slow_intro_diff}. It was shown in~\cite{BHP06} that the slow variable $x$ of this system of equations, the solution of~\eqref{e:xx}, converges in the limit as $\eps$ tends to $0$ to the solution of the homogenized equation
\begin{equation}\label{e:homog}
\dot{x} = {\bar{a}}(x) + \bar\sigma(x) \dot{W},
\end{equation}
with explicit formulas for the homogenized coefficients--see Section~\ref{sec:multisc} for details.
For finite dimensional fast systems, the coefficients in~\eqref{e:aver} and~\eqref{e:homog} can be calculated, in principle, in terms of appropriate long-time averages--see~\cite{PavlSt08} for details. The numerical method proposed in~\cite{Vand03} and analyzed in~\cite{ELV05}, coined under the name of the Heterogeneous Multiscale Method (HMM), relies on the numerical approximation of the
coefficients in ~\eqref{e:aver} and~\eqref{e:homog} by solving the original fine scale problem on time intervals of an intermediate time scale and use that data to evolve the slow variables using either~\eqref{e:aver} or~\eqref{e:homog}. In this paper we show how this methodology, when combined with a spectral method, can also be applied to SPDEs with multiple scales, that is,
to the systems \eqref{e:fast_slow_intro_adv}  and \eqref{e:fast_slow_intro_diff}. The aim of the present
paper is to present the algorithm
and report numerical experiments. The analysis of the proposed numerical method and the extension to more general classes of SPDEs with multiscale structure will be presented in a forthcoming paper.

The rest of the paper is organized as follows. In Section ~\ref{sec:algorithm} we present our new 
algorithm. Analytical and computational techniques for the analysis of SPDEs with multiple scales
at the heart of the multiscale algorithm are presented in ~\ref{sec:multisc}. In Section~\ref{sec:numerics} we present numerical experiments. Section~\ref{sec:conclusions} is reserved for conclusions and discussion on future work.

%\modif{I SUGGEST TO PROCEED DIRECTLY IN SECTION 2 WITH THE ALGORITHM
%AND LET SECTION 3 CONTAINS THE HOMOGENIZATION "JUSTIFICATION".
%SO PEOPLE JUST INTERESTED IN THE METHOD CAN READ SECT. 1/2 AND DIRECTLY JUMP TO 4}
%
%%%%%%%%%%%%%%%%%%%%%%%%%%%%%%%%%%%%%%%%%%%%%%%%%%%%%%%%%%%%%%%%%%%%%%%%%%%%%%%%%%%%%%%%
%%%%%%%%%%%%%%%%%%%%%%%%%%%%%%%%%%%%%%%%%%%%%%%%%%%%%%%%%%%%%%%%%%%%%%%%%%%%%%%%%%%%%%%%
%

\section{Numerical method}
% of the Fast/Slow System}
\label{sec:algorithm}
%\modif{
We propose a numerical algorithm to 
approximate numerically the solution of \eqref{e:spde} based on a micro-macro
algorithm, capable of capturing the effective behavior of the SPDE.
We explain the numerical algorithm for the case of diffusive time
scale (the hardest numerically) and comment on the advective time scale later in this section.
%We thus consider (see \eqref{e:spde_intro_diff})
%\begin{equation}\label{e:spde_diff}
%\partial_t u = \frac{1}{\eps^2} \cA u + \frac{1}{\eps} F(u)+ \frac{1}{\eps} Q\xi,
%\end{equation}
%where $F(u)=f(u)+\eps g(u)$.
%}
\subsection{Multiscale Algorithm}
%\modif{The following paragraph has been taken from the "theory section"}
We consider SPDEs \eqref {e:spde} 
%\begin{equation}\label{e:spde}
%\partial_t v = \cA v + f(v) + \eps Q \xi,
%\end{equation}
in a Hilbert space $\cH$ with norm $\|\cdot \|$ and inner product $\langle \cdot , \cdot\rangle$. $\cA$ denotes a differential operator, $\xi$ space-time white noise and $Q$ the covariance operator of the noise. We assume that $\cA$ is a self-adjoint nonpositive operator on $\cH$ with compact resolvent. We denote its eigenvalues and (normalized) eigenfunctions by $\{-\lambda_k, \, e_k \}_{k=1}^{\infty}$: 
\begin{equation}
\label{eig_A}
- \cA e_k = \lambda_k e_k, \quad k=1, \dots 
\end{equation}
The eigenfunctions of $\cA$ form an orthonormal basis in $\cH$. We assume that $\cA$ and the covariance operator of the noise $Q$ commute. Thus, we can write, formally, 
\begin{equation}\label{e:noise}
Q\xi = \sum_{k=1}^{+\infty} q_k e_k \xi_k(t),
\end{equation}
where $\{\xi_k(t) \}_{k=1}^{+\infty}$ are independent one-dimensional white noise processes, i.e., mean-zero Gaussian processes with $\langle  \xi_k(t)  \xi_j (s) \rangle = \delta_{k j} \delta (t-s), \, k,j=1,2,\dots$. Here
$\delta_{k j}$ and $\delta (t-s)$ are the usual Kronecker delta functions.

Furthermore, we will assume that $\cA$ has a finite dimensional kernel $\cN :=\big\{ h \in \cH \, : \, \cA h = 0  \big\}$, $\mbox{dim} (\cN) = N < + \infty$ and write $\cH = \cN \oplus \cN^{\bot}.$ We introduce the projection operators 
\begin{subequations}
\begin{eqnarray}
\label{proj_op}
\cP_c &:&\cH \mapsto \cN\\
\label{proj_op_c}
\cP_s = I - \cP_c &:& \cH \mapsto \cN^{\bot},
\end{eqnarray} 
\end{subequations}
and write $x :=P_c u$, $y := P_s u$. Finally, we will assume that noise acts only on $\cN^{\bot}$, i.e. $q_k = 0, \; k=1 \dots N$.

\noindent{\bf Step 1. Decomposition in a fast-slow system.}\\
\noindent
Using the projection operators defined in \eqref{proj_op} and \eqref{proj_op_c}, 
equation~\eqref{e:spde} can be written as a fast-slow stochastic system
\begin{subequations}\label{e:x-y}
\begin{eqnarray}\label{e:x}
\dot{x} &=& \frac{1}{\eps} \cP_c F(u), \\
\label{e:y}
\dot{y} &=& \frac{1}{\eps^2} \cA y + \frac{1}{\eps} \cP_c F(u) + \frac{1}{\eps} Q  \xi,
\end{eqnarray}
\end{subequations}
where $x(t)\in\mathbb{R}^N$ since $\mbox{dim} (\cN) = N.$ We order the pairs of eigenfuctions and eigenvalues such that the kernel $\cN$ is spanned by the first $N$ eigenfunctions of $\cA$. We can write
$$
x = \sum_{k=1}^N x_k e_k \quad \mbox{and} \quad y = \sum_{k=N+1}^{+ \infty} y_k e_k.
$$
Moreover, we introduce
\begin{eqnarray}
\label{a_proj}
a^k(x,y)&:=&\langle \cP_c F , e_k \rangle~~\hbox{ for } 1\leq k\leq N,\\
\label{b_proj}
b^k(x,y)&:=&\langle \cP_s F, e_k \rangle~~\hbox{ for } k\geq N.
\end{eqnarray}
\begin{remark}
As mentioned in the introduction we will often consider the case $F(u)=f(u)+\eps^2 g(u)$.
Then the above decomposition reads
\begin{eqnarray}
\label{a_proj_l}
a^k(x,y)&:=&\langle \cP_c f , e_k \rangle+\eps\langle \cP_c g , e_k \rangle=a_0^k(x,y)+\eps a_1^k(x,y),\\
\label{b_proj_l}
b^k(x,y)&:=&\langle \cP_s f , e_k \rangle+\eps\langle \cP_s g , e_k \rangle=b_0^k(x,y)+\eps b_1^k(x,y).
\end{eqnarray}
where we notice that for a linear function $g(u)=\nu u$ we simply have  $a_1^k(x,y)=\nu x_k,b_1^k(x,y)=\nu y_k$.
\end{remark}
Then, in view of \eqref{eig_A} and \eqref{e:noise}
we can rewrite the system~\eqref{e:x-y} in the form
\begin{subequations}\label{e:x-y-1}
\begin{eqnarray}\label{e:x-1}
\dot{x}_k &=& \frac{1}{\eps} a^k(x,y), \quad k=1, \dots N, \\
\label{e:y-1}
\dot{y}_k &=& -\frac{1}{\eps^2} \lambda_k y_k + \frac{1}{\eps}  b^k(x,y)+\frac{1}{\eps} q_k \xi_k, \quad k=N+1,N+2, \dots
\end{eqnarray}
\end{subequations}
Equations~\eqref{e:x-y}, resp. ~\eqref{e:x-y-1}, are the infinite system of singularly perturbed SDEs that we want to solve numerically.
\medskip\\
\begin{comment}
 The generator of the Markov process $\{x(t), \, y(t) \}$ is of the form
$$
\cL = \frac{1}{\eps^2} \cL_0 + \frac{1}{\eps} \cL_1 + \cL_0,
$$
where $\cL_0$ is the generator of a finite-dimensional OU processes:
$$
\cL_0 = -\sum_{j=1}^M \lambda_j y_j \frac{\partial}{\partial y_j} + \frac{1}{2} \sum_{j=1}^M q_j^1 \frac{\partial^2}{\partial y_j^2}. 
$$
Hence, the invariant measure of the fast process is Gaussian:
\begin{equation}\label{e:ou_inv_meas}
\mu(d y) = \mathcal{Z}^{-1} e^{-\sum_{j=1}^M \frac{\lambda_j y_j^2}{q_j^2}} \, dy.
\end{equation}
When the centering condition $\int f^k_{0,c}(x,y), \, \mu(dy) = 0$ is satisfied, it can be shown that for quadratic nonlinearities the slow variable $x$ converges weakly in $C([0,T]; \R^N)$ to $X$, the solution of the homogenized equation
\begin{equation}\label{e:homog}
\frac{d X}{d t} = a(X) + \sigma(X) \frac{d W}{d t},
\end{equation}
where the drift and diffusion coefficients in the homogenized equation can be calculated in terms of the solution of an appropriate Poisson equation, see~\cite[Ch. 11]{PavlSt08} and the discussion below.
\end{comment}
%
%%%%%%%%%%%%%%%%%%%%%%%%%%%%%%%%%%%%%%%%%%%%%%%%%%%%%%%%%%%%%%%%%%%%%
%
%\subsection{The Heterogeneous Multiscale Method for SDEs}
%\label{sec:HMM}
%
%
\noindent
{\bf Step 2. Truncation.}\\
\noindent
We consider a finite dimensional truncation of the above system and keep $M$ fast processes
\footnote{To simplify the notations we will use a new labeling of the index for the truncated fast system and write $(y_1,\ldots,y_M)$instead of $(y_{N+1},\ldots,y_{N+M})$ and similarly for the eigenvalues
$\lambda_k$ and the  noise intensity $q_k$.}  : 
\begin{subequations}\label{e:x-y-1_fin}
\begin{eqnarray}\label{e:x-1_fin}
\dot{\mathbf x} &=& \frac{1}{\eps} \mathbf a(\mathbf x,\mathbf y),\\
\label{e:y-1_fin}
\dot{\mathbf y} &=& -\frac{1}{\eps^2} \boldsymbol \Lambda_M\mathbf y + \frac{1}{\eps} \mathbf b(\mathbf x,\mathbf y) + \frac{1}{\eps} \boldsymbol Q_M { \boldsymbol\xi}, 
\end{eqnarray}
\end{subequations}
where ${\mathbf x}=(x_1,\ldots,x_N)^T,{\mathbf y}=(y_1,\ldots,y_M)^T,{\boldsymbol\xi}=(\xi_1,\ldots,\xi_M)^T$ and
\begin{eqnarray}
\label{a_proj_fin}
\mathbf a(\mathbf x,\mathbf y)& = &
(a^1(\mathbf x,\mathbf y),\ldots,a^N(\mathbf x,\mathbf y))^T,\\
\label{b_proj_fin}
\mathbf b(\mathbf x,\mathbf y)&=& (b^1(\mathbf x,\mathbf y),\ldots,b^M(\mathbf x,\mathbf y))^T,
\end{eqnarray}
and $\boldsymbol\Lambda_M=\hbox{diag}(\lambda_1,\ldots,\lambda_M)$ and  $\boldsymbol Q_M=\hbox{diag}(q_1,\ldots,q_M)$.
For the decomposition \eqref{b_proj_l},\eqref{a_proj_l}, we will use the notations

\begin{eqnarray}
\label{a_proj_fin_l}
\mathbf a(\mathbf x,\mathbf y)&=&\mathbf a_0(\mathbf x,\mathbf y)+\eps\mathbf a_1(\mathbf x,\mathbf y), \\
\label{b_proj_fin_l}
\mathbf b(\mathbf x,\mathbf y)&=&\mathbf b_0(\mathbf x,\mathbf y)+\eps\mathbf b_1(\mathbf x,\mathbf y),
\end{eqnarray}
where $\mathbf a_0,\mathbf a_1 \in\mathbb{R}^N$ and  $\mathbf b_0,\mathbf b_1 \in\mathbb{R}^M$ with components similar as in \eqref{a_proj_fin} or \eqref{b_proj_fin}.
\medskip\\
%
%%%%%%%%%%%%%%%%%%%%%%%%%%%%%%%%%%%%%%%%%%%%%%%%%%%%%%%%%%%%%%%
%
\noindent{\bf Step 3. Numerical solution of the reduced system.}\\
\noindent

The reduced system \eqref{e:x-y-1_fin} is solved by a micro-macro algorithm
following \cite{Vand03,ELV05}.
It consists of a macrosolver (we use the notation $X_n := X(t_n)$), chosen here to be the Euler-Maruyama scheme
\begin{equation}
\label{equ:macrosolver}
X_{n+1} =X_n+\Delta t {\bf \bar a}_M^n+ {\bar{\boldsymbol\sigma}}_M^n\Delta W_n,
\end{equation}
where $\Delta W_n$ (the Wiener increment) is $\mathcal{N}(0,\Delta t)$.
Notice that $\Delta t$ represents here a macrotime step, i.e., $\Delta t$ can be chosen much larger than $\eps$. The drift function $\bar{\mathbf a}_M^n\simeq \bar{\mathbf  a}_M(X_n)$ and diffusion  function $\bar{\boldsymbol\sigma}_M^n\simeq\bar{\boldsymbol\sigma}_M(X_n)$
appearing in  \eqref{equ:macrosolver},
recovered from a time-ensemble average, are given by
\begin{subequations}
\begin{eqnarray}
{\bf \bar a}_M^n & = & \frac{1}{K L} \sum_{j=1}^{K} \sum_{\ell=\ell_T}^{\ell_T + L-1} \partial_y {\bf a} (X_n,Y_{n,\ell,j}^{1}) Y_{n,\ell,j}^{2} \nonumber \\ \nonumber
&& + \frac{1}{K \, L}\frac{\delta t}{\eps^2} \sum_{j=1}^{K} \sum_{\ell=\ell_T}^{n_T +L-1} \sum_{\ell'=0}^{L'} \partial_x {\bf a} (X_n, Y_{n,\ell+\ell',j}^1) {\bf a}(X_n, Y_{n,\ell,j}^{1}), \\ \label{equ:an}  
\\ \nonumber
 {\bar{\boldsymbol\sigma}}_M^n (\bar{\boldsymbol\sigma}_M^n)^T & = &  \frac{1}{K \, L}\frac{2 \delta t}{\eps^2} \sum_{j=1}^{K} \sum_{\ell=\ell_T}^{\ell_T +L-1} \sum_{\ell'=0}^{L'}  {\bf a} (X_n,  Y_{n,\ell+\ell',j}^{1}) \otimes  {\bf a}(X_n,Y_{n,\ell,j}^{1}), \\ \label{equ:bn}
\end{eqnarray}
\end{subequations}
where ${Y}^{1},{Y}^{2}$ are the solutions of a suitable auxiliary system (given in \eqref{e:auxiliary} below) involving the
fast problem  \eqref{e:y-1_fin}. Here $K$ denotes the number of samples taken for the numerical calculation, $L,\, L'$ the number
of micro timesteps and $\ell_T$ a number of initial micro timesteps that are omitted in the averaging processes
to reduce the effect of transients (see below).

\medskip

\noindent{\bf Auxiliary system.}
As observed in \cite{Vand03}, for diffusive timescales, computing effective coefficients via time-averaging (relying on ergodicity), may require to solve \eqref{e:y-1_fin} over time $T={\cal O}(\varepsilon^{-2})$. To overcome this
problem, it was suggested again in \cite{Vand03} to replace the fast process in \eqref{e:y-1_fin} by
($\mathbf {\mathbf y}\simeq {\mathbf y}^{1}+\varepsilon {\mathbf y}^{2}$)
\begin{subequations}\label{e:auxiliary}
\begin{eqnarray}
\label{e:auxiliary_1}
\dot{\mathbf y}^{1} &=& -\frac{1}{\eps^2} \Lambda_M\mathbf y^{1} + \frac{1}{\eps} Q_M { \boldsymbol\xi},\\
\label{e:auxiliary_2}
\dot{\mathbf y}^{2} &=& -\frac{1}{\eps^2} \Lambda_M\mathbf y^{2} + 
\frac{1}{\eps^2} \mathbf b(\mathbf x,\mathbf y^{1}).
% +\frac{1}{\eps} \mathbf b_1(\mathbf x^N,\mathbf y^{1})
\end{eqnarray}
\end{subequations}
% with invariant measure $\nu_x^{\eps}(d y^M d z^M)$. 
%We introduce the notation
%$$
%\mathbf a^{\eps}(x,y) := \mathbf a_0 (x,y) + \eps \mathbf a_1 (x,y). 
%$$
The numerical approximations $Y^1,Y^2$ of \eqref{e:auxiliary_1} and \eqref{e:auxiliary_2}, respectively, are the functions appearing in the averaging procedure to recover the macroscopic
drift and diffusion functions (see \eqref{equ:an})-\eqref{equ:bn}).
Notice that we fix the slow variables in the system \eqref{e:auxiliary_2} at the current 
macro state $X_n$.
We use again  the Euler-Maruyama method and compute $Y^1,Y^2$ as
\begin{subequations}\label{e:microsolver}
\begin{eqnarray}
\label{e:microsolver_1}
Y_{n,\ell+1}^{1} &=& Y_{n,\ell}^{1}-\frac{\delta t}{\eps^2} \Lambda_M Y_{n,\ell}^{1} + \frac{\sqrt{\delta t}}{\eps} Q_M {\mathbf J_n},\\
\label{e:microsolver_2}
Y_{n,\ell+1}^{2} &=&  Y_{n,\ell}^{2} -\frac{\delta t}{\eps^2} \Lambda_M Y_{n,\ell}^{2} + 
\frac{\delta t}{\eps^2} \mathbf b(X_n,Y_{n,\ell}^{1}),
% +\frac{1}{\eps} \mathbf b_1(\mathbf x^N,\mathbf y^{1})
\end{eqnarray}
where ${\mathbf J_n}=\hbox{diag}(J_{n}^1,\ldots,J_{n}^M)$ and $J_{n}^k$ is a ${\cal N}(0,1)$ random variable.
The index $n$ refers to the macrotime, $t_n$.
\end{subequations}
We compute \eqref{e:microsolver_1}  over $L+L'$ microtime steps,  \eqref{e:microsolver_2}
over $L$  microtime steps to compute the  time-ensemble average \eqref{equ:an}. Notice that for the microsolver, the timestep $\delta t$ resolves the fine scale $\eps^2$. The initial values for the micro solver are taken to be (for $n\geq 1$
$$
 Y_{n,0}^{1}=Y_{n-1,\ell_T +L+L' -1 }^{1},\quad
 Y_{n,0}^{2}= Y_{n-1,\ell_T +L-1 }^{2},
$$ 
and $ Y_{0,0}^{1}= Y_{0,0}^{2}=0$ for $n=0$. 
The motivation for computing the above time averages is given in the next section.

\begin{remark}
We notice that the auxiliary system \eqref{e:auxiliary} is degenerate, since the noise in~\eqref{e:y-1_fin} is additive.\footnote{Indeed, the auxiliary system in~\cite{Vand03, ELV05} will always be degenerate, whenever the noise in the fast/slow system of SDEs that we want to solve is additive.} This implies that the results presented in~\cite[App. B]{ELV05} are not applicable in this case and a more elaborate analysis is required for proving geometric ergodicity. This analysis, based on the ergodic theory for hypoelliptic diffusions~\cite{MatSt02}, will be presented elsewhere. In the present paper we will assume that the auxiliary process~\eqref{e:auxiliary} is ergodic.
\end{remark}

\paragraph{Advective time scale.} 
%\modif{I modified the description here (why speaking about Fourier space here
%as this is the algorithm part ?)}

A similar algorithm can be derived for the advective time scale. We consider the fast-slow system~\eqref{e:fast_slow_intro_adv} that after projection and truncation reads
\begin{subequations}\label{e:x-y-1_fin_ad}
\begin{eqnarray}\label{e:x-1_fin_ad}
\dot{\mathbf x} &=& \mathbf a(\mathbf x,\mathbf y),\\
\label{e:y-1_fin_ad}
\dot{\mathbf y} &=& -\frac{1}{\eps}\boldsymbol \Lambda_M\mathbf y +\mathbf b(\mathbf x,\mathbf y) + \frac{1}{\eps} \boldsymbol Q_M { \boldsymbol\xi}, 
\end{eqnarray}
\end{subequations}
similarly as \eqref{e:x-y-1_fin}.
The macrosolver, chosen to be the Euler explicit method, is given by
\begin{equation*}
X_{n+1} =X_n+\Delta t {\bf a}_M^n,
\end{equation*}
where  the effective force ${\bf a}_M$ is given by the time average 
\begin{eqnarray}
\label{equ:an_adv}
{\bf \bar a}_M^n & = & \frac{1}{K L} \sum_{j=1}^{K} \sum_{\ell=\ell_T}^{\ell_T + L-1} {\bf a} (X_n,Y_{n,\ell,j}),\\ 
\end{eqnarray}
where $Y_{n,\ell,j}$ is a numerical approximation of the truncated fast system \eqref{fast_adv} with a slow variable
fixed at time $t_n$. As previously, $K$ denotes the number of samples and $L$ the number
of micro timesteps and $\ell_T$ is the number of initial micro timestep ommited to reduce the 
transient effects.  For the advective scaling, there is no need for an auxiliary problem
for the micro solver \cite{Vand03}.

\section{Averaging and Homogenization for SPDEs}
\label{sec:multisc}

In this section we summarize recent results on the averaging and homogenization for SPDEs~\cite{CerraiFreidlin09, BHP06} that are the analytical foundation on which the
numerical algorithm presented in Section~\ref{sec:algorithm} is built.

\subsection{Analytic form of the homogenized coefficients}
%\modif{part appearing at various places are merged here}
\label{hom}
In this section we briefly discuss the analytical form of the effective
system corresponding to \eqref{e:x-y-1_fin}.
Under the assumption that the vector field $\mathbf a_0(\mathbf x,\mathbf y)$ (see \eqref{a_proj_fin_l}) is centered with respect to the invariant measure of the fast process,
\begin{equation}\label{e:centering}
\int_{\R^M} \mathbf a_0(\mathbf x ,\mathbf y) \, \mu (d \mathbf y) = 0,
\end{equation}
then the slow process converges to a homogenized equation of the form
\begin{equation}\label{e:homog_M}
d X = \bar{\boldsymbol a}_M(X) \, dt + \bar{\boldsymbol\sigma}_M(X) \, d W,
\end{equation}
where $W$ represent an $N-$dimensional Wiener process and the SDE~\eqref{e:homog_M} is interpreted in the It\^{o} sense.
The subscript $M$ are used to emphasise the fact that the homogenized coefficients depend on the number of fast processes that we take into account.
An analytic expression for the coefficients that appear in ~\eqref{e:homog_M} is given by
%\cite{ELV05}
\begin{subequations}\label{e:coeff_homog}
\begin{eqnarray}
\bar{\mathbf a}_M( \mathbf{ x})
& = & 
\lim_{\eps \rightarrow 0} \int_{\R^M \times \R^M}  \nu_x^{\eps}(d {\mathbf y}^{1}, d {\mathbf y}^{2}) \nabla_y \mathbf a(\mathbf x,{\mathbf y}^{1}) {\mathbf y}^{2} \\ & & 
+ \lim_{\eps \rightarrow 0} \int_{\R^M} \mu (d \mathbf y^{1}) \int_0^{+\infty} \E_{y^1} \nabla_x \mathbf a(\mathbf x, 
\mathbf y_{\eps^2 s}^{1})  \mathbf a(\mathbf x,\mathbf y^{1}) \, ds, \nonumber \\
\bar{\boldsymbol \sigma}_M(\mathbf{ x})(\bar{\boldsymbol \sigma}_M(\mathbf{ x}))^T  & = & 2 \lim_{\eps \rightarrow 0} \int_{\R^M} \mu (d\mathbf y^{1}) \mathbf a(\mathbf x,\mathbf y^{1})
\nonumber \\& &
\otimes \int_0^{+\infty} \E_{y^{1}} \mathbf a(\mathbf x, \mathbf y_{\eps^2 s}^{1})ds.
%\nonumber
\end{eqnarray}
\end{subequations}
Here $\mu(d \mathbf y_{1})$ denotes the invariant measure of the process $\mathbf y^{1}$ which is given by~\eqref{e:ou_inv_meas} and $\nu_x^{\eps}(d \mathbf y^{1}, d \mathbf y^{2})$ denotes the invariant measure of the the process 
$\{\mathbf y^{1}, \, \mathbf y^{2} \}$. Notice that $\mathbf y_{\eps^2 s}^{1}=\mathbf {\tilde y}_{\tau}^{1}$ is the solution
of the rescaled process corresponding to \eqref{e:auxiliary_1}, i.e., 
$\dot{\tilde {\mathbf y}}^{1} = -\Lambda_M\mathbf {\tilde y}^{1} + Q_M { \boldsymbol\xi}$.
Alternatively, 
the calculation of the coefficients $\mathbf a_M(x)$ and $\boldsymbol \sigma_M(x)$ which appear in the homogenized equation can be obtained by the solution of the Poisson equation
\begin{equation}\label{e:poisson}
-\cL_M \phi = a_0(x, y),
\end{equation}
where $\cL_M$ is the generator of the fast (truncated) Ornstein-Uhlenbeck process. This process is ergodic and its invariant measure is Gaussian: 
\begin{equation}\label{e:ou_inv_meas}
\mu (d\mathbf y) = \frac{1}{\mathcal{Z_M}} e^{-\sum_{j=1}^M \frac{\lambda_j y_j^2}{q_j^2}} \, d\mathbf y,
\end{equation}
where $\mathcal{Z_M}$ denotes the normalization constant. We notice that the system~\eqref{e:x-y-1_fin} is
a finite dimensional fast-slow system of SDEs for which standard homogenization theory applies~\cite{lions,PapStrVar77,PavlSt08}. 
 For quadratic nonlinearities the Poisson equation~\eqref{e:poisson} can be solved analytically. The calculation of the coefficients in the homogenized (amplitude) equation reduces then to the calculation of Gaussian integrals that can also be done analytically. This will be done in Section~\ref{sec:diff_time_scale}.
%
%
%%%%%%%%%%%%%%%%%%%%%%%%%%%%%%%%%%%%%%%%%%%%%%%%%%%%%%%%%%%%%%%%%%%%%%%%%%%%%%%%%
%
%
\subsection{The Advective Time Scale}
%\modif{removed the first paragraph of this section which was redundant}
Averaging problems for fast-slow systems of SPDEs were studied recently in~\cite{CerraiFreidlin09} and their results can be applied to~\eqref{e:fast_slow_intro_adv}. One important observation is that in the system~\eqref{e:fast_slow_intro_adv}, the fast process is, to leading order $\cO(1/\eps)$, an infinite dimensional Ornstein-Uhlenbeck process. The ergodic properties of such an infinite dimensional process can be analyzed in a quite straightforward way and the invariant measure, if it exists, is a Gaussian measure in an appropriate Hilbert space that can be written down explicitly~\cite{DapZab92, DaPratoZabczyk96}.\footnote{The analysis presented in~\cite{CerraiFreidlin09} also applies to the case where the fast process is given by a semilinear parabolic SPDE. In this more general case, however, it is not possible to obtain an explicit formula for the invariant measure of the fast process.} Assuming that the process
$$                     
\partial_t z = \cA z + Q \xi
$$
is ergodic with Gaussian invariant measure $\mu$ with mean $0$ and covariance operator $\frac{1}{2}\cA^{-1}Q^2$, then the slow process $x$ converges to the solution of the averaged equation
\begin{equation}\label{e:averaged}
\dot{x} = \bar{a}(x), \quad \bar{a}(x) = \int a(x,y) \, \mu(dy),
\end{equation}
where the integration is over an appropriate Hilbert space.

When $F(\cdot)$ in~\eqref{e:spde} is given in terms of a symmetric bilinear map, i.e., $F(v)=B(v,v)$ the calculation of the vector field that appears in the averaged equation reduces to the calculation of Gaussian integrals and can be performed explicitly. In this case we have
$$P_c B(x,y):= a(x,y)= D(x,x) + C(x,y) + E(y,y),$$ where 
\begin{eqnarray*}
D_m(x,x) &=& \sum_{k, \ell =1}^N B_{k \ell m} x_k x_{\ell},\\  
C_m(x,y) &=& 2 \sum_{k=1}^N \sum_{\ell=N+1}^{\infty} B_{k \ell m} x_k y_{\ell},\\
E_m(x,y) &=& \sum_{k,\ell=N+1}^{\infty} B_{k \ell m} y_k y_{\ell}, \, m=1, \dots N,
\end{eqnarray*}
and where we used the notation $B_{k \ell m}:= \langle B(e_k, e_{\ell}), e_m \rangle$ and $N:=\mbox{dim}(\cN)$ denotes the dimension of the null space of $\cA$. 
Then, the fast-slow system~\eqref{e:fast_slow_intro_adv} becomes
\begin{subequations}\label{e:averag_quadratic}
\begin{eqnarray}
\dot{x} & = &  D(x,x) + C(x,y) + E(y,y), \\
\label{fast_adv}
\dot{y}  & = & \frac{1}{\eps} \cA y + b(x,y) + \frac{1}{\sqrt{\eps}} Q\xi,
\end{eqnarray}
\end{subequations}
and the averaged equation for~\eqref{e:averag_quadratic} reads
\begin{equation}\label{e:averaged_quadratic}
\dot{x} = D(x,x) + E,
\end{equation}
where 
$$
E_m = \sum_{k=N+1}^{+\infty} \frac{q_k^2}{2 \lambda_k} B_{k k m}, \quad m = 1, \dots N.
$$
In the case when the null space is one-dimensional, $N=1$, the averaged equation becomes
\begin{equation}\label{e:averaged_quad}
\frac{d X}{d t} = D X^2 + E, 
\end{equation}
with $D = B_{1 1 1}$ and $E_m = \sum_{k=N+1}^{+\infty} \frac{q_k^2}{2 \lambda_k} B_{k k 1} $. This equation can be solved in closed form:
$$
x(t) = \sqrt{\frac{E}{D}} \tan \left(\sqrt{E D} t + \arctan \left(\frac{D x_0}{\sqrt{E D}} \right) \right).
$$
We remark that solutions to~\eqref{e:averaged_quad}, depending on the choice of the initial conditions, do not necessarily exist for all times. We also remark that it is straightforward to consider the case where there is an additional higher order linear term (in $\eps$) in the equation, i.e. $F(v)  = B(v,v) + \eps \nu v$. In this case the averaged equation~\eqref{e:averaged_quadratic} becomes
$$
\dot{x} = D(x,x) + \nu x + E,
$$ 
where $x \in \R^N$.
%
%
%%%%%%%%%%%%%%%%%%%%%%%%%%%%%%%%%%%%%%%%%%%%%%%%%%%%%%%%%%%%%%%%%%%%%%%%%%%%%%%%%%%
%
\subsection{The Diffusive Time Scale}
\label{sec:diff_time_scale}
%\modif{I skiped the first part redundant with the introduction}
We consider the system \eqref{e:spde_intro_diff} obtained after a diffusive time rescaling
to~\eqref{e:spde}.
In order to describe the homogenized equation, we further
assume that $F(v)$ in \eqref{e:spde} is of the form 
\begin{equation}\label{e:FF}
F(v)=B(v,v)+\eps^2 \nu v,
\end{equation}
where $B(\cdot,\cdot)$ is a symmetric bilinear map satisfying $P_c B(e_k, e_k) = 0$. \footnote{This is essentially the centering condition from homogenization theory, see Equation~\eqref{e:centering} below.}

We recall that the noise does not act directly on the slow variables, $\langle Q e_k , e_k \rangle = 0, \, k=1\dots N$, where $N$ is the dimension of the null space of $\cA$.  Under appropriate  assumptions on the quadratic nonlinearity and on the covariance operator of the noise, together with the assumptions on $\cA$ and $Q$ stated earlier in this section, it is possible to prove~\cite{BHP06} that the projection of the solution to~\eqref{e:spde_intro_diff}
onto the null space of $\cA$, $x:=P_c u$, converges weakly to the solution of the homogenized SDE (the amplitude equation)
\begin{equation}\label{e:homog_multi}
d X = \bar{\mathbf a}(X) \, dt + \bar{\boldsymbol\sigma}(X) \, d W(t), \quad X(0) = X_0.
\end{equation}
where the noise is interpreted in the It\^{o} sense and the drift $  \bar{\mathbf a}(x)$ given by
\begin{equation}\label{e:F}
\bar{\mathbf a}(x) = A_\infty x - B_\infty (x,x,x) + \nu x
\end{equation}
where the linear map $A_{\infty}: \cN \rightarrow \cN$ and the trilinear map $B_{\infty}: \cN^3 \rightarrow \cN$ are defined by
\begin{subequations}
\begin{eqnarray}
A_\infty x &=&  2 B_c \Big((I \otimes_s \cA)^{-1} (B_s \otimes_s I) + (I \otimes \cA^{-1} B_s ) \nonumber 
           \\ &&+ 2 (B_c \otimes \cA^{-1})) \Big) (x \otimes \widehat{Q}), 
\label{e:nu} \\
B_\infty &=& -2  B_c (x, cA^{-1} B_s (x,x)). \label{e:eta}
\end{eqnarray}
\end{subequations}
In the above we used the notation $B_s := P_s B$ and $B_c := P_c B$, whereas $\otimes_s$ stands for the symmetric tensor product\footnote{Given a Hilbert space
$\HH$ we denote by $\HH \otimes_s \HH$ its symmetric tensor product. 
Similarly, we use the notation $v_1 \otimes_s v_2 = \frac12 \bigl(v_1\otimes v_2 + 
v_2\otimes v_1\bigr)$ for the symmetric tensor product of two elements and 
$(A \otimes_s B)(x\otimes y) = \frac12 \bigl(Ax\otimes By + By\otimes Ax\bigr)$ for the 
symmetric tensor product of two linear operators. Furthermore, we extend the bilinear form $B$ to the tensor product space 
by $B( u\otimes v)=B(u,v)$. More details can be found in~\cite[Sec. 4]{BHP06}.} and where we have defined
$$
\widehat{Q} = \sum_{k=N+1}^{\infty} \frac{q_k^2}{2 \lambda_k} \big(e_k \otimes e_k \big).
$$
The quadratic form associated with the diffusion matrix $\bar{\boldsymbol\sigma}^2$ is given by
\begin{equation}\label{e:Sigma}
\langle y, \bar{\boldsymbol\sigma}^2(x) y \rangle = 4 \sum_{k=N+1}^{+\infty} q_k^2 \langle y, B_c (e_k, x) \rangle^2 + \sum_{k, \ell=N+1}^{+\infty} \frac{q_k^2 q_{\ell}^2}{2 \lambda_{\ell} (\lambda_k + \lambda_{\ell})} \langle y, B_c (e_k, e_{\ell}) \rangle^2.
\end{equation}
Furthermore, the fast process can be approximated by an infinite dimensional Ornstein-Uhlenbeck process. The precise statement and proof of the above results can be found in~\cite{BHP06}.

\begin{remark}
The assumption that the $\cO(\eps^2)$ term in~\eqref{e:FF} is linear is needed in order to go from~\eqref{e:spde} to~\eqref{e:spde_intro_diff} after rescaling or, equivalently, to~\eqref{e:fast_slow_intro_diff}. If our starting point is the already rescaled SPDE~\eqref{e:spde_intro_diff}, then we can apply the results from~\cite{BHP06} to nonlinearities of the form $F(v) = B(v,v) + \eps^2 h(v)$ where $h(\cdot)$ is an arbitrary nonlinearity. In this case the drift term in the amplitude equation~\eqref{e:F} becomes 
\begin{equation}\label{e:Fgen}
\bar{\mathbf a}(x) = A_\infty x - B_\infty (x,x,x) + \int P_c h(x,y) \, \mu(dy),
\end{equation}
where $\mu(dy)$ denotes the invariant measure of the fast OU process.
\end{remark}

When the null space of $\cA$ is one dimensional and, consequently, the homogenized SDE is a scalar equation, it is possible to obtain sharp error estimates and to prove convergence in the strong topology. In this case Equation~\eqref{e:homog_multi} becomes

\begin{equation}\label{e:amplitude_gen}
dX=\bar a(X)+\bar\sigma(X){dW}, \quad X(0) = \langle u_0,e_1 \rangle\;,
%\dot{X} = A_{\infty} X - \frac{1}{4 \lambda_2} X^3 + \sqrt{2 \left( \frac{q_2^2}{8 \lambda_2^2} X^2 + C_{\infty}) \right) }\dot{W}.
\end{equation}
where 
\begin{equation}\label{e:amplitude_gen_coeff}
\bar a(X)=A_{\infty} X - B_{\infty}X^3,\quad \bar\sigma(X)= \sqrt{C_{\infty} + D_\infty X^2}.
\end{equation}
In the one dimensional case the formulas for the coefficients that appear in the homogenized equation have a simpler form than in the multidimensional case. In particular, we have, with $B_{k \ell m} = \langle B(e_k, e_{\ell}), e_m \rangle$:
\begin{subequations}\label{e:coeffs}
\begin{eqnarray}
A_{\infty} &= & \nu + \sum_{k=2}^\infty  {2 B_{k11}^2 q_k^2 \over \lambda_k^2}
    +  \sum_{k,\ell=2}^\infty\frac{B_{k11} B_{\ell \ell k} q_\ell^2}{\lambda_k \lambda_\ell}
        +  \sum_{k,\ell=2}^\infty 
                \frac{2B_{k\ell 1}B_{k1\ell}}{ \lambda_k+\lambda_\ell}
                \frac{q_k^2}{ \lambda_k} \;, \qquad\label{e:coeff1}\\
B_{\infty} & = &  -\sum_{k=2}^{\infty} \frac{2 B_{k11} B_{11k}}{\lambda_k}\;,\label{e:coeff2}\\
C_{\infty}&=& \sum_{m,k=2}^{\infty} \frac{2 B_{km1}^2 q_k^2 q_m^2}{(\lambda_k+\lambda_m)^2\lambda_k}\;,\quad
D_\infty =  \sum_{k=2}^\infty\frac{4 B_{k11}^2 q_k^2}{\lambda_k^2}.
\label{e:coeff3}
\end{eqnarray}
\end{subequations}
It is worth mentioning that if we are using a non-orthonormal basis, i.e. a basis $\hat{e}_k = c_k e_k$, then the coefficients that appear on the right hand side of the above equation transform according to
\begin{equation}\label{e:transformation}
\hat{B}_{k \ell m} = \frac{c_k c_{\ell}}{c_m} B_{k \ell m}.
\end{equation}
We also have $\hat{q}_k = c_k q_k$.

\begin{remark}
The formulas for the coefficients that appear in the amplitude equation~\eqref{e:homog_multi} can be also obtained by writing the SPDE~\eqref{e:fast_slow_intro_diff} in Fourier space, truncating and then using singular perturbation theory-type of techniques for the corresponding backward Kolmogorov equation~\cite{Kur73, Pap76}. More details on this approach can be found in~\cite{MTV01}. We also remark that, in general, both additive as well as multiplicative noise will appear in the amplitude equation, although only (degenerate) additive noise is present on the SPDE~\eqref{e:spde}.  
\end{remark}
%
%%%%%%%%%%%%%%%%%%%%%%%%%%%%%%%%%%%%%%%%%%%%%%%%%%%%%%%%%%%%%%%%%%%%%%%%%%%%%%%%%%%
%%%%%%%%%%%%%%%%%%%%%%%%%%%%%%%%%%%%%%%%%%%%%%%%%%%%%%%%%%%%%%%%%%%%%%%%%%%%%%%%%%%
%
%{\bf DO WE WANT TO SAY SOMETHING HERE FOR THE ADVECTIVE TIMESCALE ???}

%
%%%%%%%%%%%%%%%%%%%%%%%%%%%%%%%%%%%%%%%%%%%%%%%%%%%%%%%%%%%%%%%%%%%%%%%%%%%%%%%%%%%%%%%%%%%%%%%%%%%%%%%%
%

%
%%%%%%%%%%%%%%%%%%%%%%%%%%%%%%%%%%%%%%%%%%%%%%%%%%%%%%%%%%%%%%%%%%%%%%%%%%%%%%%%%%%%%%%%%%%%%
%
%
\section{Numerical Experiments}
\label{sec:numerics}
%\modif{some changes in  paragraphs order}
In this section we apply our numerical method to several SPDEs
and report its convergence and performance. 
We consider here several examples of SPDEs with quadratic nonlinearities
and check that the theory developed in~\cite{BHP06} and summarized in Section~\ref{sec:multisc} applies. For all of these examples we can derive rigorously the homogenized equation, with explicit formulas for the coefficients and therefore, we can present a rigorous numerical study
for our algorithm and test the  effectiveness of the proposed numerical algorithm.
%These examples thus provide us with a natural testbed for testing the effectiveness of the proposed numerical algorithm.

%As already mentioned, these examples provide us with a natural testbed for testing the effectiveness of the proposed numerical algorithm.
%In order to present a rigorous numerical study, we choose examples for which we can compute explicitly the coefficients in the amplitude equation. 
%In this way, we can compare between with the results of the numerical computations.

%First we present some analytical results that are necessary for the comparison between the theoretical results and the results of the numerical experiments.

\subsection{Theoretical considerations}
\label{sec:examples}

We will consider variants of the Burgers and the Kuramoto-Shivashinsky (KS) equations (with a linear instability term added) in one dimension with either Dirichlet or periodic boundary conditions. In particular, we will consider the singularly perturbed SPDEs (i.e. we have already rescaled to the diffusive time scale)
\begin{equation}\label{e:burgers}
\partial_t u= \frac{1}{\eps^{2}} (\partial_x^2 + 1)u +\frac{1}{\eps} u \partial_x u + \nu u 
+\frac{1}{\eps } Q \xi\;
\end{equation}
and
\begin{equation}
\label{e:KS}
\partial_t u = \frac{1}{\eps^2}(-\partial_x^2 - \partial_x^4) u + \frac{1}{\eps} u \partial_x u  + \nu u + \frac{1}{\eps} Q \xi,
\end{equation}
respectively, where the noise $\xi$ is as in Section~\ref{sec:multisc}. The operator $Q$, the covariance operator of the noise, has eigenvalues $\{ q_k \}_{k=1}^{\infty}$ and eigenfunctions $\{e_k \}_{k=1}^{\infty}$, which are also the eigenfunctions of the  differential operator that appears in either~\eqref{e:burgers} or~\eqref{e:KS}, i.e. the two operators commute. We will consider these two equations either on $[0,\pi]$ with Dirichlet boundary conditions or on $[-\pi, \pi]$ with periodic boundary conditions.
\begin{remark}
For the Burgers nonlinearity and for the boundary conditions that 
we consider it is straightforward to check that the centering
condition $P_c B(e_k, e_k) = 0$ is satisfied.
A more natural equation to consider than \eqref{e:KS} would be the KS equation in the small viscosity regime, i.e. $$\partial_t u = \frac{1}{\eps^2}(- \partial_x^2 - \mu \partial_x^4) u + \frac{1}{\eps} u \partial_x u  + \frac{1}{\eps} Q \xi, $$ where $\mu =1 - \nu, \, \nu \in (0,1)$. This equation can be rewritten in the form 
\begin{equation}\label{e:KS_small_visc}
\partial_t u = \frac{1}{\eps^2}(- \partial_x^2 -  \partial_x^4) u + \frac{1}{\eps} u \partial_x u  + \nu \partial_x^4 u + \frac{1}{\eps} Q \xi. 
\end{equation}
The theory presented in~\cite{BHP06} and the numerical scheme developed in this paper apply to this equation. The application of the numerical method developed in this paper to Equation~\eqref{e:KS_small_visc} and to related models will be presented elsewhere. Some recent analytical and numerical results on the behaviour of solutions to~\eqref{e:KS_small_visc} have been reported in~\cite{Pavliotis_al2010}.
\end{remark}
We will use the notation 
$$
\cA_B=( \partial_x^2 +1 ) \quad \mbox{and} \quad \cA_{KS} = -\partial_x^2 - \partial_x^4.
$$   
It is possible to check that for the above equations and for the chose boundary conditions the theory developed in~\cite{BHP06} and summarized in Section~\ref{sec:multisc} applies. 
Consider first equations~\eqref{e:burgers} and~\eqref{e:KS} on $[0,\pi]$ with Dirichlet boundary conditions. In this case the null space of $\cA_B$ and $\cA_{KS}$ is one dimensional: 
$$
\cN (\cA_*) = \mbox{span} \big\{\sin(\cdot) \big\}. 
$$
with $*$ being either $B$ or $KS$. The normalized eigenfunctions of $\cA_B$ and $\cA_{KS}$ are $e_k = \sqrt{\frac{2}{\pi}} \sin (\pi k)$. The corresponding eigenvalues are
$$
\lambda_k^B = k^2 -1 \quad \mbox{and} \quad \lambda_k^{KS} = k^4 - k^2,\hbox{ for $k=1,2,\ldots$.}
$$
Since the null space is one-dimensional, the homogenized equation is a scalar SDE. For the nonlinearity  $B[u,v] = \frac{1}{2} \partial_x (u v) $ it is straightforward to calculate $B_{k \ell m} = \langle B(e_k, e_{\ell}), e_m \rangle$. We have
\begin{equation}
\label{e:repBklm}
B_{k \ell m} 
= \frac{1}{2\sqrt{2\pi}} \big( |k+ \ell| \delta_{k + \ell, m} - |k - \ell| \delta_{|k -\ell|,m} \big)\;,
\end{equation}
where $\delta_{k\ell}$ denotes the Kronecker delta. We can then use formulas~\eqref{e:coeffs} to calculate the formulas that appear in the homogenized equation. Let $\{-\lambda_k \}_{k=1}^{+\infty}$ of either $\cA_{B}$ or $\cA_{KS}$ with Dirichlet boundary conditions. The homogenized equation is given by  \eqref{e:amplitude_gen}
that we recall here for convenience
\begin{equation}\label{e:amplitude_gen_bis}
dX=\bar a(X)+\bar\sigma(X){dW},
%\dot{X} = A_{\infty} X - \frac{1}{4 \lambda_2} X^3 + \sqrt{2 \left( \frac{q_2^2}{8 \lambda_2^2} X^2 + C_{\infty}) \right) }\dot{W}.
\end{equation}
%\begin{equation}\label{e:amplitude_gen}
%dX=A(X)+B(X){dW},
%%\dot{X} = A_{\infty} X - \frac{1}{4 \lambda_2} X^3 + \sqrt{2 \left( \frac{q_2^2}{8 \lambda_2^2} X^2 + C_{\infty}) \right) }\dot{W}.
%\end{equation}
where 
\begin{equation}\label{e:amplitude_gen_coeff1}
\bar a(X)=A_{\infty} X - \frac{1}{4 \lambda_2} X^3,\quad \bar\sigma(X)= \sqrt{2\left( \frac{q_2^2}{8 \lambda_2^2} X^2 + C_{\infty}) \right) }.
\end{equation}
The coefficients that appear in~\eqref{e:amplitude_gen_coeff1} can be computed as
\footnote{We use the non-normalized basis $\hat{e}_k =  \sin (\pi k)$ and use formula~\eqref{e:transformation}. }
%, when we keep infinitely many fast variables, are\footnote{We use the non-normalized basis $\hat{e}_k =  \sin (\pi \cdot)$ and use formula~\eqref{e:transformation}. }
\begin{subequations}\label{e:coeffs_homog_burgers}
\begin{eqnarray}
A_{\infty}  &=& \left( \nu +\frac{1}{8}\frac{q_2^2}{\lambda_2^2} + \frac{1}{8}\sum^{+\infty}_{k=2}\frac{k \lambda_k q_{k+1}^2 -\lambda_{k+1} q_k^2 (k+1)}{(\lambda_{k+1} + \lambda_k) \lambda_k \lambda_{k+1}}\right),
\label{e:a_inf} \\
C_{\infty} & = & \left( \frac{1}{16} \sum_{k=2}^{+\infty}\frac{q_k^2 q_{k+1}^2}{\lambda_k \lambda_{k+1}(\lambda_k + \lambda_{k+1})}\right).
\label{e:c_inf}
\end{eqnarray}
\end{subequations}
In the case where only the second mode is forced with noise, $q_2 = \sigma, \, q_M = 0, \, M=3,\dots$ then the coefficients become
$$
A_{\infty} = \nu +\frac{1}{8}\frac{\sigma^2}{\lambda_2^2} -\frac{3}{8} \frac{\sigma^2}{\lambda_2 (\lambda_2 +\lambda_3)}, \quad C_{\infty} = 0. 
$$
In this case only multiplicative noise appears in the homogenized equation and it can lead to intermittent behavior of solutions as well as noise induced transitions~\cite{Pavliotis_al2010}.

We will also consider either the Burgers or the KS equation on $[-\pi, \pi]$ with periodic boundary conditions. In this case the null space of both $\cA_B$ and $\cA_{KS}$ is two-dimensional and is spanned by 
$$
\cN(\cA_{*}) = \mbox{span} \left\{\sin(\cdot), \, \cos( \cdot) \right\},
$$
with $*$ being either $B$ or $KS$. The homogenized equation is given by  \eqref{e:homog_multi},
where ${X} = (X_1, \, X_2)$. It consists of a system  of two coupled SDEs.
We can use formulas~\eqref{e:F} and~\eqref{e:Sigma}, together with the formula for the nonlinearity $B[u,v] = \frac{1}{2} \partial_x (u \, v)$ to calculate the coefficients that appear in the homogenized equation.

\subsection{Numerical Experiments}

We shall now apply our numerical algorithm to the model problems \eqref{e:burgers}, \eqref{e:KS} described in Section \ref{sec:examples}. As the behavior of our algorithm is similar for  the Burgers and the Kuramoto-Shivashinsky equation
we will do a thorough numerical study on the Burgers equation and comment on the results for the  Kuramoto-Shivashinsky equation.
%\footnote{In fact, the stronger dissipativity properties of the Kuramoto-Shivashinsky equation %render the implementation of our numerical scheme easier.}
 %
%%%%%%%%%%%%%%%%%%%%%%%%%%%%%%%%%%%%%%%%%%%%%%%%%%%%%%%%%%%%%%%%%%%%%%%%%

\noindent{\bf Burgers Equation.}
We consider equation \eqref{e:burgers} on $[0,\pi]$ with homogeneous Dirichlet boundary conditions.
We know from Section~\ref{sec:examples} that, for $\eps$ sufficiently small, we have that
\begin{equation}
\label{approx_spde}
u(\cdot, t) \approx X(t) \sin(\cdot),
\end{equation}
where $X(t)$ is the solution of \eqref{e:amplitude_gen_bis}.
The function $\bar a(X),\, \bar\sigma(X)$ in \eqref{e:amplitude_gen_coeff} 
depends on $A_{\infty},\, C_{\infty}$ which for the Burgers equation 
can be computed using formulas~\eqref{e:coeffs_homog_burgers} with $\lambda_k = k^2 -1$. They read  $A_\infty=0.0026744369,~C_\infty=0.00026592835.$

Following the algorithm described in Section \eqref{sec:algorithm}, we look for a solution to~\eqref{e:burgers} of the form
\begin{equation}
\label{equ:n_mode_trunc}
u(\cdot,t) \simeq x(t) \sin(\cdot) + \sum_{k=1}^M y_{k}(t) \sin(k \cdot),
\end{equation}
substitute the expansion in~\eqref{e:burgers} to obtain
the a fast-slow system of SDEs as described in \eqref{e:x-y-1_fin}. Following the algorithm of Section \ref{sec:algorithm}
we compute numerically the slow variable $X_n$
as
\begin{equation}\label{e:amplitude_gen_ref_3}
X_{n+1} =X_n+\Delta t {\bar a}_M^n+ {\bar \sigma}_M^n\Delta W_n,
\end{equation}
where ${a}_M^{n}, \, {\sigma}_M^{n}$ are given by \eqref{equ:an} and \eqref{equ:bn}, respectively. We also consider the  truncated homogenized problem, i.e.,
\begin{equation}\label{e:amplitude_gen_ref_n}
d{X} = \bar a_{M}(X)dt+\bar\sigma_{M}(X)d{W},
\end{equation}
where 
\begin{equation}\label{e:amplitude_gen_coeff_bis}
\bar a_M(X)=A_M X - \frac{1}{12} X^3,\quad \bar\sigma(X)= \sqrt{2\left( \frac{1}{72} X^2 + C_M) \right) },
\end{equation}
and where $A_M, \, C_M,$ are obtained from  \eqref{e:a_inf},\eqref{e:c_inf} with the sums truncated at $M$.

For numerical comparison we also compute 
\begin{eqnarray}
\label{e:amplitude_gen_ref_4}
X_{n+1,\hbox{\scriptsize inf}}&=&X_{n,\hbox{\scriptsize inf}}+\Delta t {\bar a}(X_{n,\hbox{\scriptsize inf}})+ {\bar \sigma}(X_{n,\hbox{\scriptsize inf}})\Delta W_n,\\
\label{e:amplitude_gen_ref_5}
X_{n+1,\hbox{\scriptsize hom}}&=&X_{n,\hbox{\scriptsize hom}}+\Delta t {\bar a_M}(X_{n,\hbox{\scriptsize hom}})+ {\bar \sigma_M}(X_{n,\hbox{\scriptsize hom}})\Delta W_n,
\end{eqnarray}
the Euler-Maruyama approximation of the SDEs \eqref{e:amplitude_gen_bis} and \eqref{e:amplitude_gen_ref_n},
respectively. The same Brownian path will be used in~\eqref{e:amplitude_gen_bis}, \eqref{e:amplitude_gen_ref_3}
 and \eqref{e:amplitude_gen_ref_n}. We emphasize
that the numerical solutions for \eqref{e:amplitude_gen_bis} and \eqref{e:amplitude_gen_ref_n} rely
on analytically computed homogenized coefficients, whereas for \eqref{e:amplitude_gen_ref_3}
we implement the multiscale algorithm of Section \ref{sec:algorithm}, where the coefficients
${a}_M^{n}, {\sigma}_M^{n}$ are computed "on the fly" and rely on the microsolver
\eqref{e:microsolver_1} and \eqref{e:microsolver_2}. Hence no a-priori analytical knowledge of
the amplitude equation is required.
\begin{comment}
\noindent{\bf Version 1}
According to the results of  \cite{ELV05} for the error we have
\begin{eqnarray}
& & \mbox{sup}_{j\cdot\Delta\leq T}\left(\mathbb{E}|a_{M}^n-a_{M}(X_n)|+
\mathbb{E}|\sigma_{M}^n-\sigma_{M}(X_n)|\right),\\
& & \leq C\left(\frac{\delta t}{\eps^2}+e^{-\beta L' (\delta t/\eps^2)}\right)+
C\left(\frac{e^{-\beta n_T(\delta t/\eps^2)}}{\sqrt{L(\delta t/\eps^2)+1}}+\frac{1}
{\sqrt{KL(\delta t/\eps^2)+1}}\right),
\end{eqnarray}
where $\delta t$ is the micro time step for the micro solver, $K,L,L',n_T$ the parameters
in \eqref{equ:an}, \eqref{equ:bn} and $\beta$ is a constant entering in a dissipative condition for
the fast process \cite{ELV05} {\bf GREG: I DON'T KNOW HOW TO MAKE SENSE OF THE
$\beta,$ ANYWAY THIS RESULTS DOES NOT HOLD STRAIGHTFORWARDLY IN OUR
CASE .... WHAT SHALL WE DO ???} .
Following \cite{ELV05}, we choose $M=1$ and set the various parameters in order
to bound the error by $2^{-p}$ for various $p$. This gives $\delta t/\eps^2={\cal O}(2^{-p}),
n_T={\cal O}(1),N={\cal O}(2^{3p}),N'={\cal O}(p\cdot2^p)$ and we choose $n_T=16,N=10\cdot(2^{3p}),N'=p\cdot 2^p$. \\
\end{comment}

We choose the values of the various parameters entering in the averaging process for the computation
of $\bar a_M^n, \, \bar \sigma_M^n$ as suggested in \cite{ELV05}, i.e.,
$K=1$, $\delta t/\eps^2={\cal O}(2^{-p}),n_T={\cal O}(1),L={\cal O}(2^{3p}),L'={\cal O}(p\cdot2^p)$. 
According to  \cite{ELV05}, this guarantees (for the case of non-degenerate fast processes) that the 
error is bounded by $2^{-p}$. In our case with a degenerate fast process an error bound
is still to be established. Here we monitor such convergence numerically.
More precisely, we set $n_T=16,L=2^{3p},L'=p\cdot 2^p$
and monitor the error using
\begin{subequations}
\begin{eqnarray}
\label{num_error_finite}
E_p^M  &=& \frac{1}{N}\sum_{n=1}^N\left( |\bar a_M^n-\bar a_{M}(X_{n,\hbox{\scriptsize hom}})|+|\bar \sigma_M^{n}-\bar \sigma_{M}(X_{n,\hbox{\scriptsize hom}})|\right)
\\
E_{l,p}^M &=&\frac{1}{N}\sum_{n=1}^N\left( |\bar a_M^n-\bar a(X_{n,\hbox{\scriptsize inf}})|+|\bar \sigma_M^n-\bar\sigma(X_{n,\hbox{\scriptsize hom}})|\right),
\label{num_error_infinite}
\end{eqnarray}
\end{subequations}
for various values of $p$,
where $\Delta t=T/N$ and  $T$ represent the final time. Notice that \eqref{num_error_finite} captures the error between~\eqref{e:amplitude_gen_ref_n}--the homogenized solution of the truncated system--and the numerical solution of the truncated system, 
while~\eqref{num_error_infinite}, where the index $l$ stands for limit, captures the error 
between the homogenized solution of the limit problem~\eqref{e:amplitude_gen_bis} and the numerical solution of the truncated system.

\medskip

\noindent{\bf 2-mode truncation.} We set $M=2$ in \eqref{equ:n_mode_trunc} and substitute the expansion in~\eqref{e:burgers} to obtain
the following system of equations
\begin{subequations}\label{e:3mode}
\begin{eqnarray}
\dot{x} &=& \nu x - {1\over {2\eps}}\bigl(x y_1 + y_1 y_2 \bigr),\\
\dot{y}_1 &=& \left(\nu -\frac{3}{\eps^{2}} \right) y_1 - \frac{1}{\eps} 
              \bigl(x y_2- {1\over 2}x^2\bigr) + \frac{q_1}{\eps} 
              \xi_1(t),  \\
\dot{y}_3 &=& \left(\nu -\frac{8}{\eps^{2}} \right) y_2 + 
              {3\over {2\eps}}  \bigl(x y_1)
               + \frac{q_2}{\eps} \xi_2(t).
  \end{eqnarray}
\end{subequations}
The auxiliary process can be derived as explained in Section  \ref{sec:algorithm} and reads
\begin{subequations}\label{e:auxiliary3}
\begin{eqnarray}
\label{e:auxiliary3_11}
\dot{y}_1^{1} &=& -\frac{3}{\eps^{2}}y_1^{1}+\frac{q_1}{\eps} 
              \xi_1(t),\\
              \label{e:auxiliary3_12}
\dot{y}_2^1 &=& -\frac{8}{\eps^{2}}y_2^{1}+\frac{q_2}{\eps} 
              \xi_2(t),\\
\label{e:auxiliary3_21}
\dot{y}_1^2 &=& -\frac{3}{\eps^{2}}y_1^{2}-\frac{1}{\eps^2}
\left(x y_2^{1}-\frac{x^2}{2}\right),\\
\label{e:auxiliary3_22}
\dot{y}_2^2 &=& -\frac{8}{\eps^{2}}y_2^{1}+\frac{3}{2\eps^2}
 x y_1^{1}.
\end{eqnarray}
\end{subequations}
We apply the algorithm of Section \ref{sec:algorithm} to get a numerical approximation of the homogenised problem corresponding to \eqref{e:3mode}. The final time is $T=1$ and $N=10$, which corresponds to macro time-step of size $\Delta t=0.1$.
The macro solver for the method is given by \eqref{e:amplitude_gen_ref_3}. 
As mentioned above, we compare our results with \eqref{e:amplitude_gen_ref_4}
and \eqref{e:amplitude_gen_ref_5}.
The unknown coefficients $A_3, \, C_3$ in  \eqref{e:amplitude_gen_coeff_bis} can be computed using
\eqref{e:a_inf} and \eqref{e:c_inf}, where the sum is truncated at $M+1=3$ and read
$A_3=0.003735726834,~C_3=0.0002593873518$. 

\begin{figure}[h]
\begin{center}
\includegraphics[width=6cm]{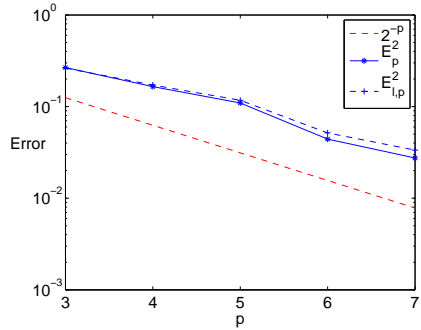}
\caption{Numerical convergence for 2-mode truncation.
On the horizontal axis we monitor the accuracy
of the micro-time step and on the vertical axis we measure the error as given by \eqref{num_error_finite} and \eqref{num_error_infinite} with $M=2$.}
\label{fig3mode}
\end{center}
\end{figure}
We observe in Figure \ref{fig3mode} that we get numerically the expected order of convergence
corresponding to $\delta t/\eps^2={\cal O}(2^{-p})$. Furthermore, as the micro time-step becomes
smaller, the numerical scheme gets closer to \eqref{e:amplitude_gen_ref_3}
and slightly deviates from \eqref{e:amplitude_gen_bis}. This is expected as the numerical
solution is not converging to that latter solution. We observe nevertheless that
with only two fast modes, the numerical scheme already captures quite well the effective
behavior of the slow variable of the infinite dimensional system.

We also illustrate the time evolution of one trajectory comparing over the time $0\leq t\leq T$
with $T=10$, the Euler-Maruyama method for the amplitude equation \eqref{e:amplitude_gen_ref_4},
the homogenized equation \eqref{e:amplitude_gen_ref_5} and the macro solver \eqref{e:amplitude_gen_ref_3}. The same Brownian path is used for generating the three trajectories and as
well as the same macro time step. 
\begin{figure}[h!]
\includegraphics[width=4.5cm,height=4cm]{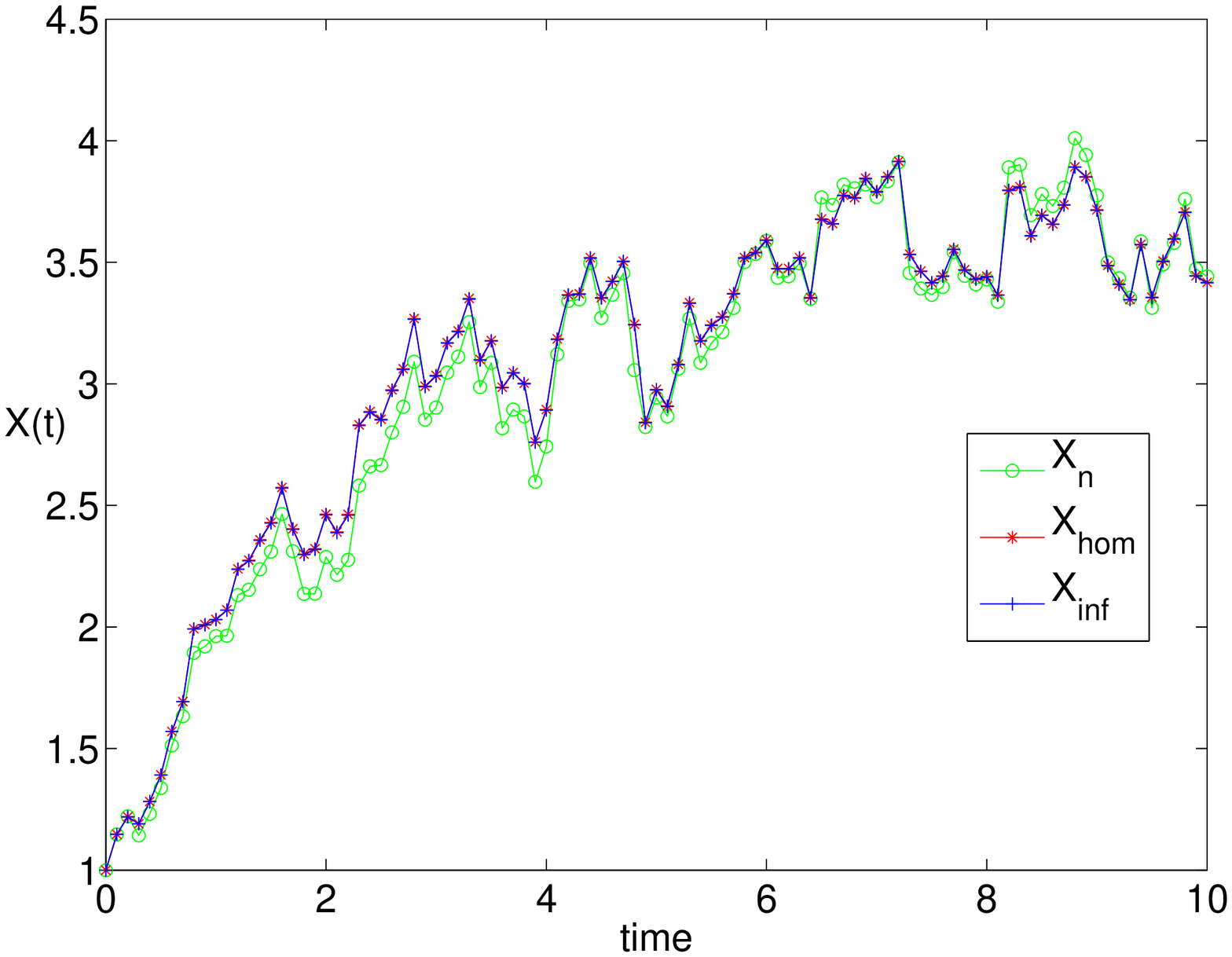}
\hskip-8mm
\includegraphics[width=4.5cm,height=4cm]{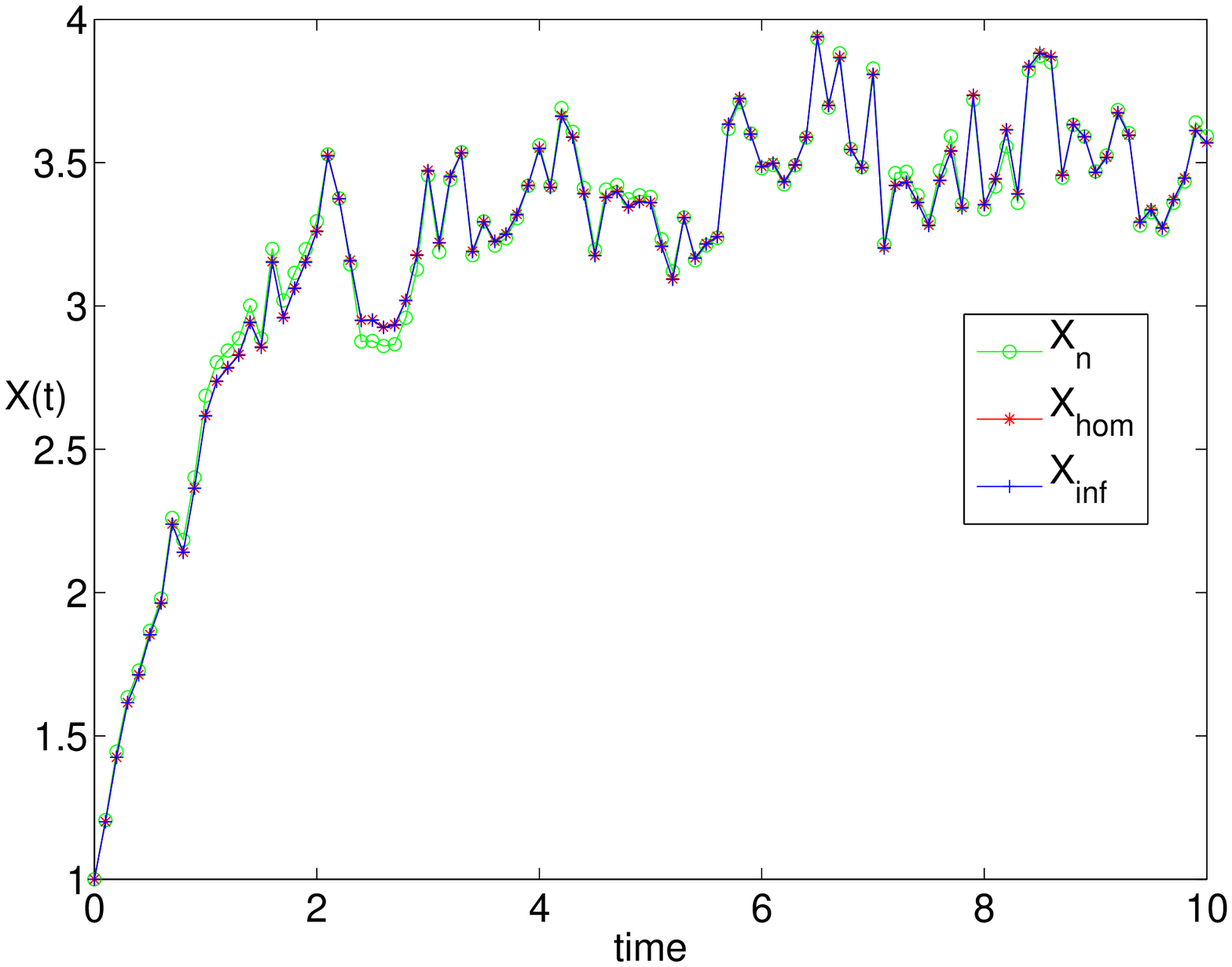}
\hskip-8mm
\includegraphics[width=4.5cm,height=4cm]{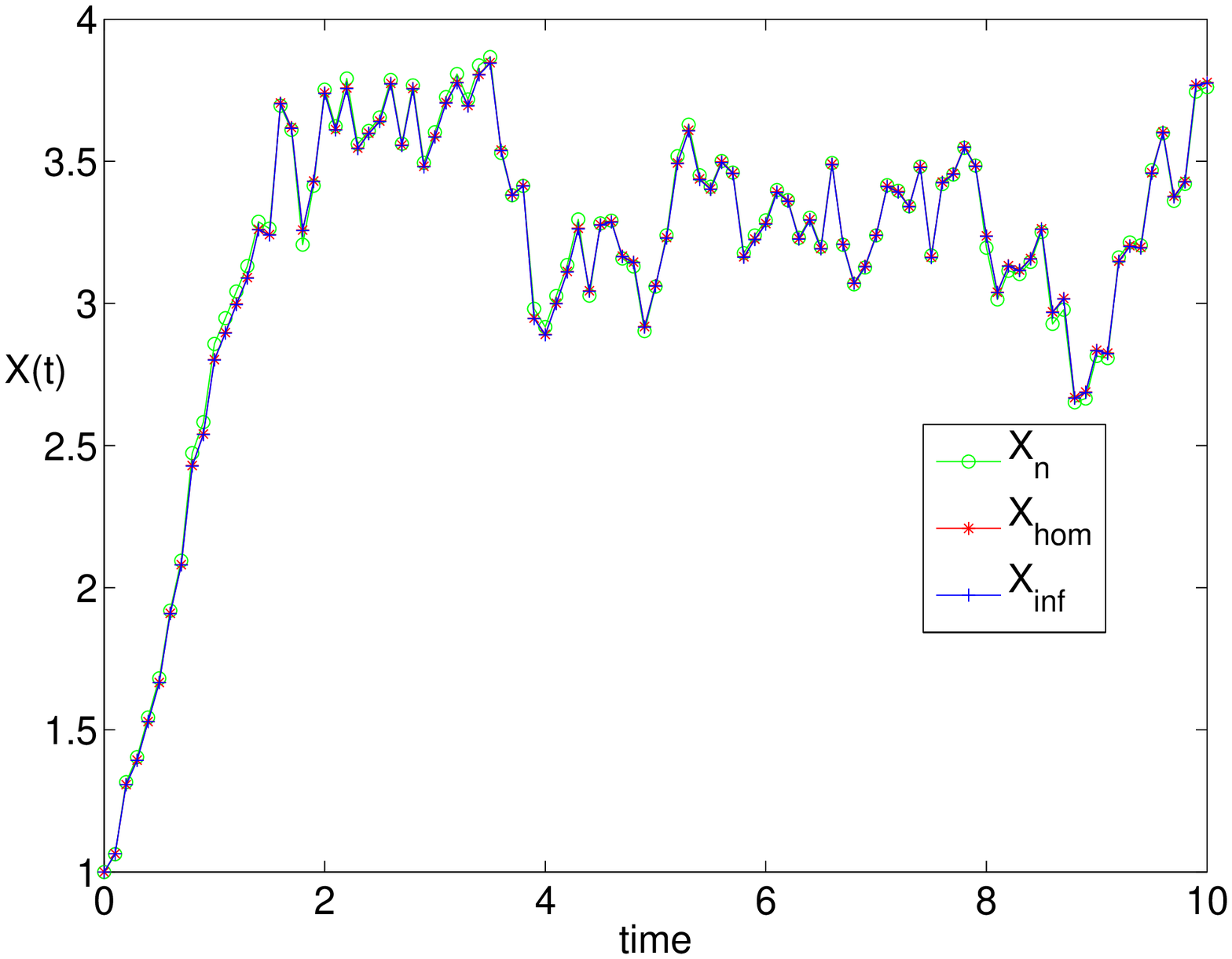}
\caption{Euler-Maruyama methods \eqref{e:amplitude_gen_ref_3} (solution denoted $X_n$),
\eqref{e:amplitude_gen_ref_5} (solution denoted $X_{n,\hbox{\scriptsize hom}}$)
and \eqref{e:amplitude_gen_ref_4} (solution denoted $X_{n,\hbox{\scriptsize inf}}$)
for three paths (left $p=3$ for $X_n$, middle $p=4$ for $X_n$, right
$p=5$ for $X_n$). We use 2-mode truncation  for \eqref{e:amplitude_gen_ref_3} and \eqref{e:amplitude_gen_ref_5}.}
\label{fig:trajectory}
\end{figure}
We perform this comparison for increasing accuracy
of the micro solver used to recover the macro data, namely, $\delta t/\eps^2={\cal O}(2^{-p}),~p=3,4,5.$
We see in Figure \ref{fig:trajectory} that the trajectory for the amplitude equation and the homogenized
equation coincide, while the macro solver gets closer to the true dynamics as we refine the micro time step.
For the same trajectory we also give a space-time plot for the approximation of the
original SPDE $
u(\cdot, t) \approx X(t) \sin(\cdot),$ with $X(t)$ solution of the  amplitude equation,
the homogenized equation or the macro solver. Again we see that the numerical
method captures the right behavior of the solution.
\begin{figure}[h!]
\begin{center}
\includegraphics[width=14cm,height=7cm]{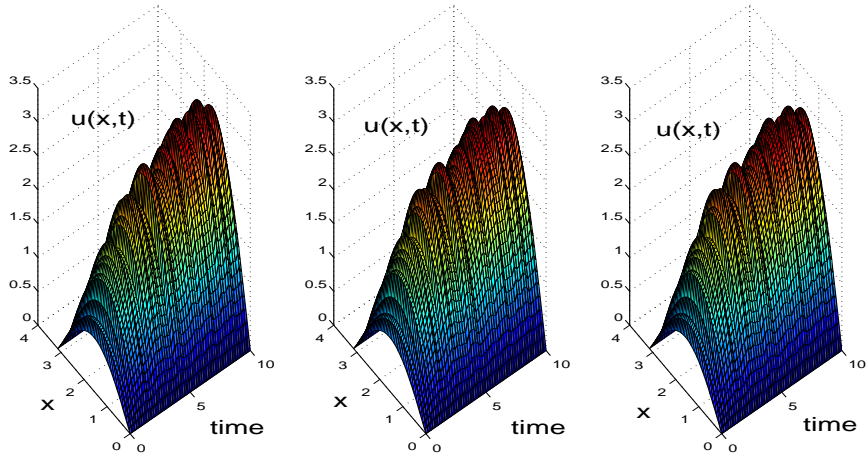}
\caption{Approximation \eqref{approx_spde} of the solution $u(x,t)$ of the SPDE;
$u(\cdot,t)\simeq X_n(t) \sin(\pi\cdot)$ (left figure $p=3$), $u(\cdot,t)\simeq X_{n,\hbox{\scriptsize hom}}(t) \sin(\pi\cdot)$
(middle figure) and 
$u(\cdot,t)\simeq X_{n,\hbox{\scriptsize inf}}(t) \sin(\pi\cdot)$ (right figure).}
\label{spacetime1}
\end{center}
\end{figure}
%\modif{I removed Figs. 4 and 5 I don't think it is necessary to have three time
%the space-time comparison, as the picture looks the same}
%\begin{figure}[h!]
%\begin{center}
%\includegraphics[width=15cm,height=8cm]{figure/3modes/time_space_3mode_t10p4.eps}
%\caption{Approximation \eqref{approx_spde} of the solution $u(x,t)$ of the SPDE;
%$u(x,t)\simeq X_n \sin(\cdot)$ (left figure $p=4$), $u(x,t)\simeq X_{n,\hbox{\scriptsize hom}}$
%(middle figure) and 
%$u(x,t)\simeq X_{n,\hbox{\scriptsize inf}}$ (right figure).
%The brownian path
%is different from the brownian path of figure \ref{spacetime1}.}
%\label{spacetime2}
%\end{center}
%\end{figure}
%\begin{figure}[h!]
%\begin{center}
%\includegraphics[width=15cm,height=8cm]{figure/3modes/time_space_3mode_t10p5.eps}
%\caption{Approximation \eqref{approx_spde} of the solution $u(x,t)$ of the SPDE;
%$u(x,t)\simeq X_n \sin(\cdot)$ (left figure $p=5$), $u(x,t)\simeq X_{n,\hbox{\scriptsize hom}}$
%(middle figure) and 
%$u(x,t)\simeq X_{n,\hbox{\scriptsize inf}}$ (right figure).
%The brownian path
%is different from the brownian path of figures \ref{spacetime1} and  \ref{spacetime2}.}
%\label{spacetime3}
%\end{center}
%\end{figure}
%\medskip\\

\medskip

\begin{figure}[h]
\begin{center}
\includegraphics[width=6cm]{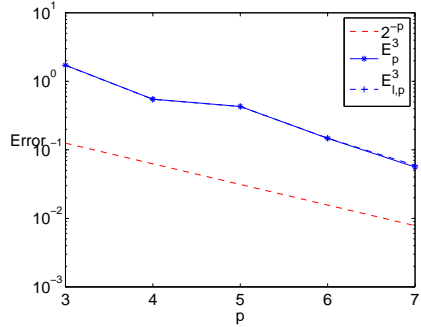}
\caption{Numerical convergence for 3-mode truncation.
On the horizontal axis we monitor the accuracy
of the micro-time step and on the vertical axis we measure the error as given by \eqref{num_error_finite} and \eqref{num_error_infinite} with $M=3$.
}
\end{center}
\label{fig4mode}
\end{figure}

\noindent{\bf 3-mode truncation.} We set $M=3$ in \eqref{equ:n_mode_trunc} and obtain
the following system of equations
\begin{subequations}\label{e:4mode}
\begin{eqnarray}
\dot{x} &=& \nu x - \frac{1}{2\eps}\bigl(x y_1 + y_1 y_2 + y_2 y_3\bigr),\\
\dot{y}_1 &=& \left(\nu -\frac{3}{\eps^{2}} \right) y_1 - \frac{1}{\eps} 
              \bigl(x y_2+y_1 y_3 - {1\over 2}x^2\bigr) + \frac{q_1}{\eps} 
              \xi_1(t),  \\
\dot{y}_2 &=& \left(\nu -\frac{8}{\eps^{2}} \right) y_2 - 
              \frac{3}{2\eps}  \bigl(x y_3 - x y_1\bigr)
               + \frac{q_2}{\eps} \xi_2(t),\\
\dot{y}_3 &=& \left( \nu - \frac{15}{ \eps^{2}} \right) y_3 + 
                    \frac{1}{\eps}\left(2 x y_2 +  y_1^2  \right)+  \frac{q_3}{ \eps } \xi_3(t).
\end{eqnarray}
\end{subequations}
The auxiliary process can be computed similarly as for the 3-mode truncation.
We perform the same set of numerical experiments as for the 3-mode truncation,
reported
in Figure \ref{fig4mode}. Similar behavior than previously noted can be observed. Observe that the discrepancy between
the numerical scheme and \eqref{e:amplitude_gen_bis} gets smaller. This is expected as with additional modes,
the homogenized equation \eqref{e:amplitude_gen_ref_n} (that we aim at approximating with our multiscale scheme)
gets closer to \eqref{e:amplitude_gen_bis}.

\medskip

\noindent{\bf 4-mode truncation.} We set $M=4$ in \eqref{equ:n_mode_trunc} and apply
the similar procedure as previously. For the sake of brevity, we do not write the system of equations in this case and just report the numerical convergence. 
%A=0.00297454122;B=1/72;C=0.0002651000164;

\begin{figure}[h]
\begin{center}
\includegraphics[width=6cm]{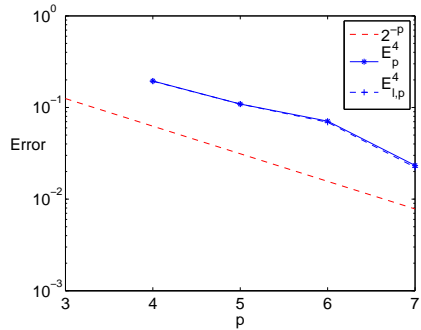}
\caption{Numerical convergence for 4-mode truncation.
On the horizontal axis we monitor the accuracy
of the micro-time step and on the vertical axis we measure the error as given by \eqref{num_error_finite} and \eqref{num_error_infinite} with $M=4$.
}
\end{center}
\label{fig5mode}
\end{figure}

We see in Figure \ref{fig5mode} a similar
behavior of our numerical scheme as observed previously. We again notice that the discrepancy between
the numerical scheme and \eqref{e:amplitude_gen_bis} is smaller than for lower order truncation.
Notice that the first numerical result reported is obtained for $\delta t/\eps^2=2^{-p}$ with $p=4$. This is due to stability issues with the Euler-Maruyama scheme for the fast process. As the linear term in the equation for the fifth mode $y_5$ is $\left( \nu - \frac{24}{ \eps^{2}}\right) y_5$, the stability restriction $24\delta t/\eps^2 \leq2$ implies $\delta t/\eps^2\leq 1/12$ and thus the time step  $\delta t/\eps^2=1/8$
is unstable.
\medskip\\

\paragraph{ The Kuramoto-Shivashinsky equation.} 
The equations for the $M$-mode truncation of the Kuramoto-Shivashinsky equation are very similar to the ones for the Burgers equation and will not be presented here. The only difference is that the fast process is more dissipative than for the Burgers equation, due to the stronger dissipativity of the operator $\cA_{KS}$ compared to $\cA_{B}$. 
% In theory the remark below is true since we converge more rapidly to the invariant measure
% but numerically there might be other issues
%This stronger dissipativity implies that our numerical method performs even better for the KS equation.
As the results of the numerical experiments for the KS equation are very similar to the results reported in this section for the Burgers equation, they will not be presented here.
\section{Conclusions and Further Work}
\label{sec:conclusions}
%\modif{Some modification done in the conclusion}
We have presented
a new numerical method for the efficient and accurate solution of stochastic partial differential equations with multiple scales. The new numerical scheme is based on a combination of a spectral method with the HMM methodology and has been tested on SPDEs with quadratic nonlinearities for which a rigorous homogenization theory exists. This enables us to check the performance of our method. The numerical experiments presented in this paper suggest that the new method performs well and allows to solve accurately multiscale SPDEs by solving a low dimensional fast-slow system of SDEs. The method is suitable for infinite dimensional stochastic systems for which there is clear separation of scales, and for which a low dimensional homogenised (or averaged) equation
for the slow modes exists.
\begin{comment}
{\bf ASSYR, CAN YOU PLEASE CHECK THE NEXT SENTENCE?} It appears that the methodology proposed in~\cite{ELV05}, combined with a spectral method, is particularly well suited for dissipative singularly perturbed SPDEs, since for such a class equations the stiffness of the resulting fast/slow systems of SDEs (due to the dissipativity and the scale separation) can be turned to our advantage.
\end{comment}

There are still many questions that are left open. First, the rigorous analysis of the proposed method and a careful study of the convergence and stability properties of the proposed method remains to be done. In addition, the optimisation of the proposed method by tuning appropriately the parameters of the algorithm has not been performed yet. This appears to be an open problem even when the HMM methodology is applied to finite dimensional fast/slow systems of SDEs~\cite{Liu2010}. 

The proposed numerical algorithm could be used to study in detail the qualitative and quantitative properties of solutions to SPDEs with quadratic nonlinearities, since SPDEs of this form exhibit very rich dynamical behaviour. Furthermore, we would like to apply the numerical algorithm to more general classes (and systems) of semilinear SPDEs, for which an averaged or homogenised equation is known to exist. Examples include systems of reaction/diffusion equations that were considered in~\cite{CerraiFreidlin09} as well as the Swift-Hohenberg SPDE~\cite{BHP04}. 

In our algorithm, we did not make use of the fact that the form of the amplitude equation (i.e. a Landau equation with additive and multiplicative noise) is known. Knowledge of the functional form of the coefficients that appear in the homogenised or averaged equation can be used in order to simplify the numerical algorithm. The stochastic Landau equation appears as the amplitude equation for several infinite dimensional stochastic dynamical systems, not only for SPDEs with quadratic nonlinearities, e.g.~\cite{BlMPSc01}. Thus, the algorithm proposed in this paper, 
could be modified to develop an efficient method for studying the dynamics of SPDEs near bifurcation. All these topics are currently under investigation.

%Third, the study of several different distinguished limits for SPDEs with quadratic nonlinearities. 
%{\bf CAN YOU GIVE A BIT MORE DETAILS} Fourth, the study of this problem {\bf THIS IS A BIT VAGUE CAN YOU GIVE SOME MORE DETAILS} on the advective time scale, since the averaged equation can exhibit finite time singularities.

% Study different distinguished limits, study numerics at the advection time scale.
% Other SPDEs for which an amplitude equation exists such as Swift-Hohenberg.
% Analysis of the method, proof of convergence, consistency etc.
% Systems of SPDEs/SPDEs in higher spatial dimensions.

\paragraph{Acknowledgments} Part of this work was done while GP was visiting the Mathematics Section of EPFL. The hospitality of the department and of the group of A. Abdulle is greatly acknowledged. The authors thank S. Krumscheid for an extremely careful reading of an earlier version of the paper and for many useful remarks. GP is supported by the EPSRC.

%
%%%%%%%%%%%%%%%%%%%%%%%%%%%%%%%%%%%%%%%%%%%%%%%%%%%%%%%%%%%%%%%%%%%%%%%%%%%%%%%%%
%\bibliography{mybib}
%\bibliographystyle{plain}
%\bibliographystyle{elsarticle-num}

\end{document}